# A Proof of the Riemann Hypothesis and Determination of the Relationship Between Non- Trivial Zeros of Zeta Functions and Prime Numbers

Murad Ahmed Abu Amr
MSc Degree in Physics, Mutah University
murad.abuamr@yahoo.com

## ABSTRACT

This analysis which uses new mathematical methods aims at proving the Riemann hypothesis and figuring out an approximate base for imaginary non-trivial zeros of zeta function at very large numbers, in order to determine the path that those numbers would take. This analysis will prove that there is a relation links the non-trivial zeros of zeta with the prime numbers, as well as approximately pointing out the shape of this relationship, which is going to be a totally valid one at numbers approaching infinity.

## INTRODUCTION

The Riemann hypothesis states that the real part of every non-trivial zeros of zeta function is always constant and equals ½. The mathematicians hope that the proof of this hypothesis to contribute in defining the locations of prime numbers on the number line.
This proof will depend on a new mathematical analysis that is different from the other well-known methods; It`s a series of mathematical steps that concludes a set of mathematical relationships resulted from conducting mathematical analyses concerning Riemann zeta function and Dirichlet eta function after generalizing its form, as I will explain hereinafter, as well as identifying new mathematical functions which I will define its base later.

I will resort to this different analysis because the known methods used in separating the arithmetic series to find its zeros cannot figure out the non-trivial zeros of zeta function; besides, I will explain the reason for its lack of success.

### *The main steps of this proof are as follows:*

(1) Identifying the mathematical method for separating the base series elements:

$$\sum_{n=1}^{N} n^{-s}, \qquad Re[s] < 1, N \to \infty$$

Where Re[s] is a real part of s.

This way of separation depends mainly on the sifting process, but we work on adding other conditions that include more ratios of numbers that have been removed by the sifting process.

(2) Adopting approximations that depend on limiting at infinity.
(3) Initially to get mathematical equations as afunction of very large prime numbers, then to generalize these equations to include all the very large natural numbers.

### *1) The ratio of numbers that have been removed through the sifting process:*

First of all, I will describe the group which includes the arithmetic sequence on which the sifting process is per formed, assuming group Ω, which represents the natural numbers, beginning with one and up to the natural number $N$ as follows:

$\Omega = \{1,2,3,4, \dots, N\}$

Assuming that $N$ is a very large prime numbers, I will start by removing the first prime numbers (2) and all its multiples from the arithmetic sequence. So, the removed sequence will be as follows:

$$\Omega_2 = \{2, 4, 6, 8, \dots, N-1\}$$

Where $\Omega_2$ is a set of even integers, so the remained group after taking the even integers will be as follows:

$$\Omega - \Omega_2 = \{1, 3, 5, 7, \dots, N\}$$

It is obvious that the numbers in each group come close to ½ $N$ whenever $N$ approaches to infinity; this matter can be mathematically formulated as follows:

$$N_2 \sim \frac{1}{2}N \Rightarrow \lim_{N \to \infty} \frac{N_2}{N} = \frac{1}{2}$$

Since $N_2$ is the number of elements in Group $\Omega_2$.

I will use the "~" sign referring to the approximation which becomes more accurate when the number becomes larger, while it reaches the total value when the number $N$ leads to infinity.

Now, I will define Group $\Omega_3$ as a set of numbers includes the second prime numbers (3) and all its remained multiples after removing $\Omega_2$ that will be as follows:

$$\Omega_3 = \{3, 9, 15, 21, \dots, j_3\}$$

Since $j_3$ is the largest number for the multiples of 3 that has been removed from the group.

Whereas the number of elements in this group is equal to one-third of what is remained in group $\Omega$ after removing $\Omega_2$, so the number of its elements equals one-third of what is left after taking ½ $N$. This can be mathematically represented as follows:

$$N_3 \sim \frac{1}{3}\left(1 - \frac{1}{2}\right)N \sim \frac{1}{6}N \Rightarrow \lim_{N \to \infty} \frac{N_3}{N} = \frac{1}{6}$$

Since $N_3$ is the number of elements in Group $\Omega_3$.

The remained elements of the parent group after taking the numbers prescribed in the last group will be as follows:

$$\Omega - \Omega_2 - \Omega_3 = \{1, 5, 7, 11, 13, \dots, N\}$$

By the same way, the number of elements in Group $\Omega_5$ - resulted from the removal of the third prime numbers (5) and its remained multiples - can be expressed after removing $\Omega_2$ and $\Omega_3$ from $\Omega$; this referred number is given the symbol $N_5$ and the number will lead to one-fifth of the remained after taking the multiples of $N_2$ & $N_3$ from $N$, as follows:

$$N_5 \sim \frac{1}{5}\left(1 - \frac{1}{2} - \frac{1}{6}\right)N \sim \frac{1}{15}N$$

In a more splendid expression:

$$N_5 \sim \frac{1}{5}\left(1 - \frac{1}{2}\right)\left(1 - \frac{1}{3}\right)N \Rightarrow \lim_{N \to \infty} \frac{N_5}{N} = \frac{1}{15}$$

By the same way, we indicate the number of elements in Group $\Omega_7$ which is $N_7$ as follows:

$$N_7 \sim \frac{1}{7}\left(1 - \frac{1}{2}\right)\left(1 - \frac{1}{3}\right)\left(1 - \frac{1}{5}\right)N \sim \frac{4}{105}N$$

If $P_n$ refers to the prime numbers $n$, where $P_1 = 2$, $P_2 = 3$, $P_3 = 5$, and so on.

Thus, the number of elements in any group $\Omega_{Pn}$ can be indicated, and number of elements in such group is symbolized by $N_{Pn}$ as follows:

$$N_{P_n} \sim \frac{1}{P_n}\left(1 - \frac{1}{2}\right)\left(1 - \frac{1}{3}\right)\left(1 - \frac{1}{5}\right)\dots\left(1 - \frac{1}{P_{n-1}}\right)N$$

Or at the following form:

$$\frac{N_{P_n}}{N} \sim \frac{1}{P_n} \prod_{k=1}^{n-1}\left(1 - \frac{1}{P_k}\right) \ldots\ldots\ldots (1)$$

After the end of the sifting process, which expires upon the removal of the last prime numbers of the group, which is $N$, only one single number will be left which is 1.

We will present a set of the first ten elements resulting from substitution in the last equation, which represents the ratios resulting from the sifting process.

$$\left\{\frac{1}{2}, \frac{1}{6}, \frac{1}{15}, \frac{4}{105}, \frac{8}{385}, \frac{16}{1001}, \frac{192}{17017}, \frac{3072}{323323}, \frac{55296}{7436429}, \frac{110592}{19605131}\right\}$$

Whereas the number $N$ equals the total number of group elements plus 1 as follows:

$$N = 1 + N_2 + N_3 + N_5 + N_7 + \cdots + N_N$$

By dividing the terms of the equation into $N$, it will take the following form:

$$1 = \frac{1}{N} + \frac{N_2}{N} + \frac{N_3}{N} + \frac{N_5}{N} + \frac{N_7}{N} + \cdots + \frac{N_N}{N}$$

Taking the limit when $N$ approaches to infinity, we can figure out this sequence which its right hand side represents the total ratio of the numbers removed due to the sifting process, and we also can prove that its summation is equals 1:

$$\frac{1}{2} + \frac{1}{3}\left(1 - \frac{1}{2}\right) + \frac{1}{5}\left(1 - \frac{1}{2}\right)\left(1 - \frac{1}{3}\right) + \frac{1}{7}\left(1 - \frac{1}{2}\right)\left(1 - \frac{1}{3}\right)\left(1 - \frac{1}{5}\right) + \cdots = 1$$

In other mathematical word:

$$\sum_{n=1}^{\infty} \frac{1}{P_n} \prod_{k=1}^{n-1}\left(1 - \frac{1}{P_k}\right) = 1$$

❖ But, What is the importance of ratios of numbers that are removed by the sifting process?

The answer is: to find out how a series of infinite numbers is being sifted by them in respect of zeta function, in case the series is divergent. Now, let`s show how the sifting process is done in zeta function and how the convergence affects on its success:

$$\zeta(t) = \sum_{n=1}^{\infty} n^{-t}, \quad Re[t] > 1 \quad \ldots\ldots\ldots (2)$$

The sifting process is done in a mathematical way as follows:

**First**: By multiplying zeta function by $2^{-t}$, then subtracting the result from the original function, as follows:

$$\zeta(t) - 2^{-t}\zeta(t) = 1 + 2^{-t} + 3^{-t} + 4^{-t} + \cdots - (2^{-t} + 4^{-t} + 6^{-t} + 8^{-t} + \cdots)$$

In this way, all factors of the prime numbers 2 are being removed from zeta function, so the equation takes the following form:

$$(1 - 2^{-t})\zeta(t) = 1 + 3^{-t} + 5^{-t} + 7^{-t} + \cdots$$

**Second**: By multiplying the tow hand side of the last equation by $3^{-t}$, then subtracting the result from the last series as in the first step; in this way, we remove all the factors of the prime numbers 3 from the series, as follows:

$$(1 - 3^{-t})(1 - 2^{-t})\zeta(t) = 1 + 5^{-t} + 7^{-t} + 11^{-t} + \cdots$$

Then we continue the same way, working on removing the factors of the prime numbers one by one until infinity, leaving only one number in the right hand side of the equation, and then it could be proved that zeta function is given by the following relation:

$$\zeta(t) = \prod_{n=1}^{\infty} \frac{1}{1 - P_n^{-t}}, \quad Re[t] > 1 \quad \ldots\ldots\ldots (3)$$

But what is the reason for the success of the sifting process in reaching this remarkable result when the series is converge?

This question can be answered through the following analysis:

Let`s go back to the first step, When we multiple the series $\zeta(t)$ by $2^{-t}$ times, all numbers involved in this series become even integers, but they reach values twice faster than the original series; so, the result of subtracting the series from itself after being multiplying by $2^{-t}$ will be at the following form:

$$(1 - 2^{-t})\zeta(t) = \lim_{k \to \infty} \left\{ \sum_{n=1}^{k}(2n-1)^{-t} + \sum_{n=1}^{k}(2n)^{-t} - \sum_{n=1}^{2k}(2n)^{-t} \right\}$$

In other formulation:

$$(1 - 2^{-t})\zeta(t) = \lim_{k \to \infty} \left\{ \sum_{n=1}^{k}(2n-1)^{-t} - \sum_{n=k+1}^{2k}(2n)^{-t} \right\}$$

In case of convergence:

$$\lim_{k \to \infty} \sum_{n=1}^{k} n^{-t} = \lim_{k \to \infty} \sum_{n=1}^{2k} n^{-t}, \quad Re[t] > 1 \quad \ldots\ldots\ldots (4)$$

Or:

$$\lim_{k \to \infty} \sum_{n=k}^{2k} n^{-t} = 0, \quad Re[t] > 1 \quad \ldots\ldots\ldots (5)$$

Therefore, the result of the first sifting process returns back to the form:

$$(1 - 2^{-t})\zeta(s) = 1 + 3^{-t} + 5^{-t} + 7^{-t} + \cdots = \sum_{n=1}^{\infty}(2n-1)^{-t}$$

While in case of divergence – when the real part of $t < 1$ - equation no. 4 is no longer valid, each sums of it has a different and divergent value, and the limit of the equation leads to infinity.

2) *The general base of zeta function:*

For the formula of zeta function to be inclusive of all values of $t$ which have a real part greater than zero, excluding whether $t = 1$, this can be done through the following method:

$$\zeta(t) = \sum_{n=1}^{\infty} n^{-t}, \quad Re[t] > 1 \Rightarrow \zeta(t) = \frac{1 - 2^{1-t}}{1 - 2^{1-t}} \sum_{n=1}^{\infty} n^{-t}, \quad Re[t] > 1$$

This formula becomes in this form after the analysis:

$$\zeta(t) = \frac{1}{1 - 2^{1-t}} \sum_{n=1}^{\infty} \frac{(-1)^{n+1}}{n^t}, \quad Re[t] > 0, t \neq 1 \quad \ldots\ldots\ldots (6)$$

Since this formula of zeta function is known between 0 and 1, so this is a more comprehensive form than the first simple form.

But how can zeta equals zero? Of course, this could be achieved when the numerator becomes equal to zero, where the numerator is equal to eta function and it can be given by the following mathematical relationship:

$$\eta_2(t) = \sum_{n=1}^{\infty} \frac{(-1)^{n+1}}{n^t}, \qquad Re[t] > 0 \qquad \ldots\ldots\ldots (7)$$

Where $\eta_2(t)$ is eta function as afunction of $t$, while the number 2 indicates that the base function here depends on the number 2, and the series which it represents can be indicated as follows:

$$\eta_2(t) = 1 - 2^{-t} + 3^{-t} - 4^{-t} + 5^{-t} - 6^{-t} + \cdots, \qquad Re[t] > 0$$

For eta function to be equal to zero, the sum of even terms must equal to the sum of odd terms, and through this the zeros of $\eta_2(t)$ can be extracted as follows:

$$1 - 2^{-t} + 3^{-t} - 4^{-t} + 5^{-t} - 6^{-t} + \cdots = 0$$

$$\Rightarrow 1 + 3^{-t} + 5^{-t} + \cdots = 2^{-t} + 4^{-t} + 6^{-t} + \cdots$$

Since the right hand side of the equation represents the result of multiplying the series $\zeta(t)$ – in its simplest form – by $2^{-t}$, the last equation becomes in the following form:

$$(1 - 2^{1-t})(1 + 2^{-t} + 3^{-t} + 4^{-t} + 5^{-t} + \cdots) = 0$$

To solve this equation, as the series is:

$$1 + 2^{-t} + 3^{-t} + 4^{-t} + 5^{-t} + \cdots \neq 0$$

Thus, the zeros will not appear unless when:

$$1 - 2^{1-t} = 0$$

- The initial solution (primary):

The initial solution shows the real part of the variable $t$, which provides a direct solution to the last equation as follows:

$$2^{1-t} = 1$$

$$\Rightarrow 1 - Re[t] = 0$$

$$\Rightarrow Re[t] = 1$$

Such value which appeared through the direct solution of the equation gives noting but the real part of the solution, and in order to extract the whole solution in both its real and imaginary parts, we have to turn the left hand side of the equation into a Sine function form, as follows:

$$1 - 2^{1-t} = 0 \Rightarrow 2^{\frac{1}{2}(1-t)} \left\{ 2^{-\frac{1}{2}(1-t)} - 2^{\frac{1}{2}(1-t)} \right\} = 0$$

$$\Rightarrow e^{\frac{1}{2}(1-t)\ln 2} - e^{-\frac{1}{2}(1-t)\ln 2} = 0$$

$$\Rightarrow \frac{1}{2i} \left\{ e^{\frac{i}{2}(1-t)\ln 2^{-i}} - e^{-\frac{i}{2}(1-t)\ln 2^{-i}} \right\} = 0, \qquad i = \sqrt{-1}$$

$$\Rightarrow \sin \left\{ \frac{1}{2}(1-t)\ln 2^{-i} \right\} = 0$$

$$\Rightarrow -\frac{i}{2}(1-t)\ln 2 = \pm j\pi, \qquad j = 1,2,3,\ldots$$

$$\Rightarrow t_j = 1 \pm \frac{2j\pi}{\ln 2} i, \qquad j = 1,2,3,\ldots \infty$$

3) **The reason for the non-appearance of non-trivial zeros in eta function through this analysis:**

The zeros that arose before us through this classical analysis are not the non-trivial zeros of zeta function, for sure, because the factor 1-2[1-t] which was the reason for the existence of those zeros –exists in the denominator of zeta function; which leads to the removal of that factor from the numerator and the disappearance of these zeros, keeping only the non-trivial zeros of zeta function.

I will call these zeros that have arisen in the last equation "the classical zeros", because the classical analysis is able to reveal them, and in order to distinguish them from the trivial zeros.

Nevertheless, the non-trivial zeros of zeta function are still the zeros of eta function as well, but we need an advanced mathematical analysis to be able to figure them out.

Returning back to the method used in extracting the zeros of zeta function, the method depends on the relation no. 5, but this relation becomes invalid when $t$ or a real part of $t$ is less than 1. Hence, in order to perform the correct analysis that can leads to a correct mathematical relation, in this case we have to take the sifting equations which we extracted into consideration.

Now I will define the following function:

$$\zeta(N) = \frac{1}{1-2^{1-s}} \sum_{n=1}^{N} \frac{(-1)^{n+1}}{n^s} \quad \ldots\ldots\ldots (8)$$

Where $s$ is an non-trivial zero of zeta function, this value $\zeta(N)$ will lead to zero whenever the number $N$ approaches to infinity.

When expanding the summation, it becomes in the following form:

$$\zeta(N) = \frac{1}{1-2^{1-s}} \{1 - 2^{-s} + 3^{-s} - 4^{-s} + 5^{-s} - 6^{-s} + \cdots N^{-s}\}$$

Now if $N$ is an even integer, the summation will definitely equal zero if the value of $s$ is zero, while it is equal to one if $N$ is an odd integer; so, zero is not considered as one of the function zeros that are defined to us.

Thus, we must redefine the function to prove that zero is not a solution in all cases, and this can be done through the following:

$$\zeta(s) = \frac{1}{1-2^{1-s}} \lim_{N \to \infty} \sum_{n=1}^{N} \frac{(-1)^{n+1}}{n^s}$$

But the zeta function can be defined as follows:

$$\zeta(s) = \frac{1}{1-2^{1-s}} \lim_{N \to \infty} \sum_{n=1}^{N+1} \frac{(-1)^{n+1}}{n^s}$$

By adding the two equations:

$$2\zeta(s) = \frac{1}{1-2^{1-s}} \lim_{N \to \infty} \sum_{n=1}^{N} \frac{(-1)^{n+1}}{n^s} + \frac{1}{1-2^{1-s}} \lim_{N \to \infty} \sum_{n=1}^{N+1} \frac{(-1)^{n+1}}{n^s}$$

$$\Rightarrow \zeta(s) = \frac{1}{2} \left\{ \frac{1}{1-2^{1-s}} \lim_{N \to \infty} \sum_{n=1}^{N} \frac{(-1)^{n+1}}{n^s} + \frac{1}{1-2^{1-s}} \lim_{N \to \infty} \sum_{n=1}^{N+1} \frac{(-1)^{n+1}}{n^s} \right\}$$

Whereas:

$$\sum_{n=1}^{N+1} \frac{(-1)^{n+1}}{n^s} = \sum_{n=1}^{N} \frac{(-1)^{n+1}}{n^s} + \frac{(-1)^N}{(N+1)^s}$$

The equation can be shortened to the following form:

$$\Rightarrow \zeta(s) = \frac{1}{1-2^{1-s}} \lim_{N \to \infty} \left( \frac{(-1)^N N^{-s}}{2} + \sum_{n=1}^{N} \frac{(-1)^{n+1}}{n^s} \right), \quad Re[s] \geq 0 \quad ....(9)$$

$N^{-s}$ was replaced with $(N+1)^{-s}$ in the first term inside the bracket of the last equation because the grow of 1 for more than a very large number such $N$ means nothing in such case; now with this definition of zeta function, the function value when $s$ = zero, equals the negative half whether $N$ is an even or an odd number, thanks to entering the factor:

$$\frac{(-1)^N N^{-s}}{2(1-2^{1-s})}$$

Now, by going back to equation no. 8 which will take the following form:

$$\zeta(N) = \frac{(-1)^N N^{-s}}{2(1-2^{1-s})} + \frac{1}{1-2^{1-s}} \sum_{n=1}^{N} \frac{(-1)^{n+1}}{n^s}$$

As previously mentioned $\zeta(N) \to 0$ when $N$ approaches to infinity, so if $N$ approaches to infinity, the sign "=" could be put instead of "~". Therefore, comes out the following equation:

$$\frac{(-1)^N}{2N^s} + \sum_{n=1}^{N} \frac{(-1)^{n+1}}{n^s} = 0$$

By neglecting the first term due to approaching to zero when $N$ approaches to infinity, the final equation will be as follows:

$$\sum_{n=1}^{N} \frac{(-1)^{n+1}}{n^s} = 0$$

$$\Rightarrow 1 + 3^{-s} + 5^{-s} + \cdots + N^{-s} = 2^{-s} + 4^{-s} + 6^{-s} + \cdots + N^{-s} \quad ......(10)$$

Of course if the left hand side of the equation ends with $N^{-s}$, the right hand side must ends with the term $(N\pm1)^{-s}$, but due to the large number of $N$ the two values are considered approximately equal. So, both hand side of the last equation end with $N^{-s}$, considering the following:

$$(N \pm 1)^{-s} = N^{-s}, \quad N \to \infty$$

We can conclude a very important thing from equation no.10, that the sum of even terms equal sum of odd terms when the number of odd terms equals the number of even terms, when $s$ is one of the non-trivial zeros of zeta function.

Therefore, the sum of both right and left hand sides equals the double of each hand side separately. This can be prescribed through the following equation:

$$1 + 2^{-s} + 3^{-s} + 4^{-s} + 5^{-s} + \cdots + N^{-s} = 2(2^{-s} + 4^{-s} + 6^{-s} + \cdots + N^{-s})$$

By taking $2^{-s}$ as a common factor from the right hand side, we get the following relation:

$$1 + 2^{-s} + 3^{-s} + 4^{-s} + \cdots + N^{-s} = 2^{1-s}\{1 + 2^{-s} + 3^{-s} + 4^{-s} + \cdots + (N/2)^{-s}\} ......(11)$$

I will define the function $\lambda(N)$ as follows:
$$\lambda(N) = 1 + 2^{-s} + 3^{-s} + 4^{-s} + 5^{-s} + \cdots + N^{-s} \quad ........(12)$$
By using this definition, equation no.11 can be formulated as follows:
$$\lambda(N) = 2^{1-s} \lambda(N/2)$$
Through the last equation we can conclude the next one:
$$\lambda(2N) = 2^{1-s} \lambda(N) ......(13)$$
By the same way as in equation no. 6, the zeta function can be defined by the following method:

$$\zeta(t) = \frac{1-3^{1-t}}{1-3^{1-t}} \sum_{n=1}^{\infty} n^{-t}, \quad Re[t] > 1$$

Through entering the factor 1-$3^{1-t}$ to the series after ordering the terms arises the following:

$$\zeta(t) = \frac{1 + 2^{-t} - 2 \times 3^{-t} + 4^{-t} + 5^{-t} - 2 \times 6^{-t} + \cdots}{1 - 3^{1-t}}, \qquad Re[t] > 0$$

Now we have a different form of eta function represents the numerator of the last equation, which is given from the following relation:

$$\eta_3(t) = 1 + 2^{-t} - 2 \times 3^{-t} + 4^{-t} + 5^{-t} - 2 \times 6^{-t} + \cdots, \qquad Re[t] > 0$$

This function $\eta_3(t)$ is the eta function that relies on the number 3 in its base, and such form of function has classical zeross arose through the next relation:

$$t_j = 1 \pm \frac{2j\pi}{\ln 3} i, \qquad j = 1,2,3,\ldots$$

We can conclude the base of the function $\eta_4(t)$ according to the pattern that brings out through the two previous bases, as follows:

$$\eta_4(t) = 1 + 2^{-t} + 3^{-t} - 3 \times 4^{-t} + 5^{-t} + 6^{-t} + 7^{-t} - 3 \times 8^{-t} + \cdots, Re[t] > 0$$

Such form of function has classical zeross arose through the next relation:

$$t_j = 1 \pm \frac{2j\pi}{\ln 4} i, \qquad j = 1,2,3,\ldots$$

In general, the general base of the function $\eta_m(t)$, which is eta function that relies at its base on the natural number $m$, can be given through the following relation:

$$\eta_m(t) = \sum_{k=1}^{\infty} \left\{ -(m-1)(km)^{-t} + \sum_{n=(k-1)m+1}^{km-1} n^{-t} \right\}, \ldots (14)$$

And the general base for the classical zeros of this function is given as follows:

$$t_j = 1 \pm \frac{2j\pi}{\ln m} i, \qquad j = 1,2,3,\ldots \qquad , \ldots\ldots\ldots (15)$$

Therefore, the zeta function can be generalized at the following form:

$$\zeta(t) = \frac{1}{1 - m^{1-t}} \sum_{k=1}^{\infty} \left\{ -(m-1)(km)^{-t} + \sum_{n=(k-1)m+1}^{km-1} n^{-t} \right\}$$

By turning back to the form of zeta function resulting from entering the factor 1-$3^{1-t}$ when $t = s$:

$$\zeta(N) = \frac{1}{1 - 3^{1-s}} \{1 + 2^{-s} - 2 \times 3^{-s} + 4^{-s} + 5^{-s} - 2 \times 6^{-s} + \cdots N^{-s}\}$$

Since $s$ is one of the non-trivial zeros of zeta function, so the function $\zeta(N)$ will lead to zero whenever the number $N$ approaches to infinity.

Whereas $s$ is one of the non-trivial zeros of zeta function, the numerator in the defined function defined to us leads to zero, as follows:

$$1 + 2^{-s} - 2 \times 3^{-s} + 4^{-s} + 5^{-s} - 2 \times 6^{-s} + \cdots N^{-s} = 0$$

$$\Rightarrow 1 + 2^{-s} + 3^{-s} - 3 \times 3^{-s} + 4^{-s} + 5^{-s} + 6^{-s} - 3 \times 6^{-s} + \cdots N^{-s} = 0$$

$$\Rightarrow 1 + 2^{-s} + 3^{-s} + 4^{-s} + 5^{-s} + \cdots N^{-s} = 3 \times 3^{-s} + 3 \times 6^{-s} + 3 \times 9^{-s} + \cdots + 3N^{-s}$$

$$\Rightarrow 1 + 2^{-s} + 3^{-s} + 4^{-s} + \cdots + N^{-s} = 3^{1-s}\{1 + 2^{-s} + 3^{-s} + 4^{-s} + \cdots + (N/3)^{-s}\}$$

Through this analysis we can obtain the following relation:

$$\lambda(N) = 3^{1-s} \lambda(N/3) \Rightarrow \lambda(3N) = 3^{1-s} \lambda(N)$$

Through this pattern, the last base can be generalized as follows:

$$\lambda(k\,N) = k^{1-s} \lambda(N)$$

Thus, it becomes clear that $\lambda(N)$ proportional to $N^{1-s}$. From this notice we can conclude the following relation:

$$\lambda(N) = a\,N^{1-s} \ldots\ldots (16)$$

Where ($a$) is a constant number that independent on $N$, the value of the constant a can be figured out through the next proof:

$$\sum_{n=1}^{N} n^{-s} = 1 + 2^{-s} + 3^{-s} + 4^{-s} + 5^{-s} + \cdots + N^{-s}$$

By ordering the terms inversely:

$$\sum_{n=1}^{N} n^{-s} = N^{-s} + (N-1)^{-s} + (N-2)^{-s} + (N-3)^{-s} + \cdots + \left(N - (N-1)\right)^{-s}$$

And by taking $N^{-s}$ as a common factor:

$$\sum_{n=1}^{N} n^{-s} = N^{-s} \left\{ 1 + \left(1 - \frac{1}{N}\right)^{-s} + \left(1 - \frac{2}{N}\right)^{-s} + \left(1 - \frac{3}{N}\right)^{-s} + \cdots + \left(1 - \frac{N-1}{N}\right)^{-s} \right\}$$

now by work on expand each of terms by using Taylor Expansion, as follows:

First term:

$$\left(1 - \frac{1}{N}\right)^{-s} = 1 + \frac{s}{N} + \frac{s(s+1)}{2!}\frac{1}{N^2} + \frac{s(s+1)(s+2)}{3!}\frac{1}{N^3} + \frac{s(s+1)(s+2)(s+3)}{4!}\frac{1}{N^4} + \cdots$$

Second term:

$$\left(1 - \frac{2}{N}\right)^{-s} = 1 + s\left(\frac{2}{N}\right) + \frac{s(s+1)}{2!}\left(\frac{2}{N}\right)^2 + \frac{s(s+1)(s+2)}{3!}\left(\frac{2}{N}\right)^3 + \frac{s(s+1)(s+2)(s+3)}{4!}\left(\frac{2}{N}\right)^4 + \cdots$$

Third term:

$$\left(1 - \frac{3}{N}\right)^{-s} = 1 + s\left(\frac{3}{N}\right) + \frac{s(s+1)}{2!}\left(\frac{3}{N}\right)^2 + \frac{s(s+1)(s+2)}{3!}\left(\frac{3}{N}\right)^3 + \frac{s(s+1)(s+2)(s+3)}{4!}\left(\frac{3}{N}\right)^4 + \cdots$$

And so on until $N$-1 of the terms.

The next step in re-ordering the terms in the series to gather every term of the terms that have similar factors in one series, as follows:

$$\sum_{n=1}^{N} n^{-s} = N^{-s} \left\{ 1 + \sum_{1}^{N-1} 1 + \frac{s}{N}\sum_{1}^{N-1} n + \frac{s(s+1)}{2!}\frac{1}{N^2}\sum_{1}^{N-1} n^2 + \frac{s(s+1)(s+2)}{3!}\frac{1}{N^3}\sum_{1}^{N-1} n^3 + \cdots \right\} \dots (17)$$

As for the series that arose in the last equation, the sum base for each of them is known as follows:

$$\sum_{n=1}^{N-1} 1 = N - 1$$

$$\sum_{n=1}^{N-1} n = \frac{1}{2}N(N-1)$$

$$\sum_{n=1}^{N-1} n^2 = \frac{1}{6}N(N-1)(2N-1)$$

$$\sum_{n=1}^{N-1} n^3 = \frac{1}{4}N^2(N-1)^2$$

And so on until the infinity of the series, but where $N$ is a very large number and approaches to infinity, the bases of the last grand series can be simplified in order to lead to the following form:

$$\sum_{n=1}^{N-1} 1 = N - 1 \Rightarrow \lim_{N \to \infty} \sum_{n=1}^{N-1} 1 = N$$

$$\sum_{n=1}^{N-1} n = \frac{1}{2}N(N-1) \Rightarrow \lim_{N \to \infty} \sum_{n=1}^{N-1} n = \frac{N^2}{2}$$

$$\sum_{n=1}^{N-1} n^2 = \frac{1}{6}N(N-1)(2N-1) \Rightarrow \lim_{N \to \infty} \sum_{n=1}^{N-1} n^2 = \frac{N^3}{3}$$

$$\sum_{n=1}^{N-1} n^3 = \frac{1}{4}N^2(N-1)^2 \Rightarrow \lim_{N \to \infty} \sum_{n=1}^{N-1} n^3 = \frac{N^4}{4}$$

Through noticing this pattern we can conclude the next relation:

$$\lim_{N \to \infty} \sum_{n=1}^{N-1} n^k = \frac{N^{k+1}}{k+1}$$

By the substitution with these endings in equation no.17, we get the following relation:

$$\sum_{n=1}^{N} n^{-s} = N^{-s} \left\{ N + \frac{s}{N} \frac{N^2}{2} + \frac{s(s+1)}{2!} \frac{1}{N^2} \frac{N^3}{3} + \frac{s(s+1)(s+2)}{3!} \frac{1}{N^3} \frac{N^4}{4} + \cdots \right\}$$

That can be simplified to the following form:

$$\sum_{n=1}^{N} n^{-s} = N^{-s} \left\{ N + s \frac{N}{2} + \frac{s(s+1)}{2!} \frac{N}{3} + \frac{s(s+1)(s+2)}{3!} \frac{N}{4} + \cdots \right\}$$

By considering $N$ as a common factor:

$$\sum_{n=1}^{N} n^{-s} = N^{1-s} \left\{ 1 + \frac{s}{2} + \frac{s(s+1)}{3!} + \frac{s(s+1)(s+2)}{4!} + \cdots \right\}$$

By dividing the right hand side into $s-1$ and multiply by the same number:

$$\sum_{n=1}^{N} n^{-s} = \frac{N^{1-s}}{s-1} \left\{ s - 1 + \frac{(s-1)s}{2} + \frac{(s-1)s(s+1)}{3!} + \frac{(s-1)s(s+1)(s+2)}{4!} + \cdots \right\}$$

The next step is to add 1 to the series inside the bracket, by then subtracting it as follows:

$$\sum_{n=1}^{N} n^{-s} = \frac{N^{1-s}}{s-1} \left\{ \left( 1 + (s-1) + \frac{(s-1)s}{2!} + \frac{(s-1)s(s+1)}{3!} + \frac{(s-1)s(s+1)(s+2)}{4!} + \cdots \right) - 1 \right\} \dots \dots (18)$$

Now I will define the following series:

$$\sigma = 1 + (s-1) + \frac{(s-1)s}{2!} + \frac{(s-1)s(s+1)}{3!} + \frac{(s-1)s(s+1)(s+2)}{4!} + \cdots$$

Thus, equation no. 18 can be written as follows:

$$\sum_{n=1}^{N} n^{-s} = \frac{N^{1-s}}{s-1} \{\sigma - 1\} \dots \dots (19)$$

To be able to determine the value of the series $\sigma$ we have to resort to the series that represents Taylor Expansion for the term $(1-x)^{1-s}$:

$$(1-x)^{1-s} = 1 + (s-1)x + \frac{(s-1)s}{2}x^2 + \frac{(s-1)s(s+1)}{3!}x^3 + \frac{(s-1)s(s+1)(s+2)}{4!}x^4 + \cdots$$

The base of this series will be typical to the series σ when $x \to 1$ as follows:

$$\sigma = \lim_{x \to 1}(1-x)^{1-s}$$

But the value of $(1-x)^{1-s}$ depends on the domain of the real part of s as follows:

$$\lim_{x \to 1}(1-x)^{1-s} = \begin{cases} 0, & Re[s] < 1 \\ \pm\infty, & Re[s] > 1 \end{cases} \quad \ldots \ldots (20)$$

So that, the value of σ which provide an converge solution is zero, and through this equation no.19 can be written as follows:

$$\sum_{n=1}^{N} n^{-s} = \frac{N^{1-s}}{s-1}(0-1) \Rightarrow \sum_{n=1}^{N} n^{-s} = \frac{N^{1-s}}{1-s}, N \to \infty, Re[s] < 1 \ldots \ldots (21)$$

Through this result, we can conclude the following results:

First: identifying the value of the constant ($a$) by comparing equation no.21 with equation no.16, and through the comparison, we can find that the constant ($a$) is given through the following relation:

$$a = \frac{1}{1-s}$$

From this constant ($a$) we can figure out the base of the series $\lambda(N)$ as follows:

$$\lambda(N) = \frac{N^{1-s}}{1-s} \quad \ldots \ldots \ldots (22)$$

Second: identifying the domain where the real part of non-trivial zeros of zeta function exist. As prescribed in equation no. 20, the domain where the zeros exist is $Re[s]<1$.

### 4) *The ratio of series resulting from the sifting numbers With prime numbers factors:*

Now I will prove that the ratio of the sums removed by sifting to the ratio of the main sum equals the ratio that I have determined in the first part of this article (equation no.1).

Returning back to equation no.10, it is shown that the sum based on odd integers starting from 1 and up to the number N is equal to the sum based on even integers starting from N and up to $N \pm 1$:

$$1 + 3^{-s} + 5^{-s} + \cdots + N^{-s} = 2^{-s} + 4^{-s} + 6^{-s} + \cdots + N^{-s}$$

Hence, the left hand side equals half of $\lambda(N)$, since that of $\lambda(N)$ equals the sum of both series together. This can be expressed as follows:

$$1 + 3^{-s} + 5^{-s} + \cdots + N^{-s} = \frac{1}{2}\lambda(N)$$

The left hand side can be divided as follows:

$$1 + 3^{-s} + 5^{-s} + \cdots + N^{-s} = \{1 + 5^{-s} + 7^{-s} + \cdots + N^{-s}\} + \{3^{-s} + 9^{-s} + 15^{-s} + \cdots + N^{-s}\}$$

But the second bracket in the right hand side in last equation is the result of sifting the prime numbers 3 with its remaining multiples after removing the multiples of 3, which in turn is divided as follows:

$$3^{-s} + 9^{-s} + 15^{-s} + \cdots + N^{-s} = 3^{-s}\{1 + 3^{-s} + 5^{-s} + \cdots + (N/3)^{-s}\} = 3^{-s} \times \frac{1}{2}\lambda(N/3)$$

$$\Rightarrow 3^{-s} + 9^{-s} + 15^{-s} + \cdots + N^{-s} = 3^{-s} \times \frac{1}{2} \times 3^{s-1} \lambda(N) = \frac{1}{6}\lambda(N)$$

Now I will define the following functions:

$$\lambda_2(N) = 2^{-s} + 4^{-s} + 6^{-s} + \cdots + j_2^{-s}$$

$$\lambda_3(N) = 3^{-s} + 9^{-s} + 15^{-s} + \cdots + j_3^{-s}$$

$$\lambda_5(N) = 5^{-s} + 25^{-s} + 35^{-s} + \cdots + j_5^{-s}$$

The number $j_2 = N -1$, while the number $j_3$ was previously defined, it is the largest number of the multiples of the prime numbers 3 that exists in the series $\lambda(N)$ and it is not one of the multiples of number 2. The same way with the number $j_5$, as it is the largest number of the multiples of the prime numbers 5 that exists in the series $\lambda(N)$ and it is not one of the multiples of the prime numbers 2 or 3, and so on.

As it is shown that:

$$\lambda_2(N) = \frac{1}{2}\lambda(N) \Rightarrow \lim_{N\to\infty} \frac{\lambda_2(N)}{\lambda(N)} = \frac{1}{2}$$

Also:

$$\lambda_3(N) = \frac{1}{6}\lambda(N) \Rightarrow \lim_{N\to\infty} \frac{\lambda_3(N)}{\lambda(N)} = \frac{1}{6} = \frac{1}{3}\left(1 - \frac{1}{2}\right)$$

By the same way we can prove that:

$$\lambda_5(N) = \frac{1}{15}\lambda(N) \Rightarrow \lim_{N\to\infty} \frac{\lambda_5(N)}{\lambda(N)} = \frac{1}{15} = \frac{1}{5}\left(1 - \frac{1}{2}\right)\left(1 - \frac{1}{3}\right)$$

Since it has been shown through this pattern that the ratio between $\lambda_{Pn}$ functions to $\lambda(N)$ equal the ratios resulting from the sifting process, so the result can be generalized as follows:

$$\lambda_{P_n}(N) = \frac{1}{P_n}\left(1 - \frac{1}{2}\right)\left(1 - \frac{1}{3}\right)\left(1 - \frac{1}{5}\right)\cdots\left(1 - \frac{1}{P_{n-1}}\right)\lambda(N)$$

$$\Rightarrow \lim_{N\to\infty} \frac{\lambda_{P_n}(N)}{\lambda(N)} = \frac{1}{P_n}\prod_{k=1}^{n-1}(1 - P_k^{-1}) \quad \ldots\ldots\ldots (23)$$

5) **Finding the non-trivial zeros:**

Through Riemann`s functional equation:

$$\zeta(t) = 2^t\, \pi^{t-1} \sin\left(\frac{\pi t}{2}\right)\Gamma(1-t)\zeta(1-t) \quad \ldots\ldots\ldots (24)$$

Through this equation the non-trivial zeros of zeta function was appear, and they represent all values of $t$ which equals a negative even integers; this is through the factor of $\zeta(1-t)$, so the non-trivial zeros of the function $\zeta(t)$ are the zeros of both $\zeta(1-t)$ and $\zeta(t)$ together.

Therefore, the non-trivial zeros of the function $\zeta(t)$ are the values of $t$ which make $\zeta(1-t)$ = zero. Whereas $\zeta(1-t)$ and $\zeta(t)$ have the same function base, the two functions become equal when $t=1-t$, including cases that the functions equal to zero.

Also if $\zeta(t)$ can be given through equation no.6, $\zeta(1-t)$ can be given through the following relation:

$$\zeta(1-t) = \frac{1}{1-2^t}\sum_{n=1}^{\infty}\frac{(-1)^{n+1}}{n^{1-t}}, \quad Re[t] < 1, t \neq 0 \quad \ldots\ldots\ldots (25)$$

It is clearly shown from examining the last two relations that the non-trivial zeros, if it found, have a real part exists during the interval at which the real domain of the two functions intersects and which is the interval lies on between zero and one, as follows:

$$\zeta(t) \cap \zeta(1-t), Re[t] \in (0,1)$$

Assuming the following:

$$t = s = x + iy$$

Where $s$ is one of the non-trivial zeros of zeta function, $x$ is the real part of $s$, while $y$ is the imaginary part of it, and $i$ is the root of -1.

Thus:

$1 - s = 1 - x - iy$

The next step is to substitute with $x$ and $y$ in $s$ at numerator only, as follows:

$$\Rightarrow \zeta(s) = \frac{1}{1 - 2^{1-s}} \sum_{n=1}^{\infty} (-1)^{n+1} n^{-x-iy} = \frac{1}{1 - 2^{1-s}} \sum_{n=1}^{\infty} (-1)^{n+1} n^{-x} \{\cos(y \ln n) - i \sin(y \ln n)\}$$

What is important here is the numerator, as it contains the zeros and it is eta function. The real part can be separated from the imaginary part of eta function as follows:

$$Re[\eta_2(s)] = \sum_{n=1}^{\infty} (-1)^{n+1} n^{-x} \cos(y \ln n)$$

$$Im[\eta_2(s)] = -\sum_{n=1}^{\infty} (-1)^{n+1} n^{-x} \sin(y \ln n)$$

Where $Im[\eta_2(t)]$ is imaginary part of $\eta_2(t)$.

As for the function $\zeta(1-t)$:

$$\zeta(1 - s) = \frac{1}{1 - 2^s} \sum_{n=1}^{\infty} (-1)^{n+1} n^{x-1+iy} = \frac{1}{1 - 2^s} \sum_{n=1}^{\infty} (-1)^{n+1} n^{x-1} \{\cos(y \ln n) + i \sin(y \ln n)\}$$

$$\Rightarrow Re[\eta_2(1 - s)] = \sum_{n=1}^{\infty} (-1)^{n+1} n^{x-1} \cos(y \ln n)$$

$$\Rightarrow Im[\eta_2(1 - s)] = \sum_{n=1}^{\infty} (-1)^{n+1} n^{x-1} \sin(y \ln n)$$

For both functions $\zeta(1-s)$ and $\zeta(s)$ equal to zero, their real part and imaginary part must equal to zero. When comparing the real parts in both functions and the imaginary parts as well, we notice that every two bases are similar in form, except for the following:

(1) In the sign concerning the imaginary part and this in turn will not affect on the equality when the part is equal to zero.
(2) In the power of number $n$, and here the two functions become equal when the following condition is achieved:

$-x = x - 1 \Rightarrow x = \frac{1}{2}$

The last proof does not prove the existence of non-trivial zeros of zeta function and that its real part equals the half, but it proves that if there are non-trivial zeros of zeta function, its real part does not take any value except the half.

Hence, if we were able to prove the existence of non-trivial zeros when performing equations on the very large numbers, the last proof proves its generalization in all cases.

❖ *Using λ functions in finding the non-trivial zeros of zeta function:*

Let's go back to $\lambda_{Pn}$, each function of it represents the result of sifting one of the prime numbers, so the function $\lambda_{Pn}$ shall take the following from:

$$\lambda_{P_n}(N) = P_n^{-s} + P_n^{-2s} + (P_n\,P_{n+1})^{-s} + \cdots + j_{P_n}^{-s} \quad \ldots\ldots\ldots (26)$$

Where $j_{Pn}$ is the largest number of the multiples of the prime numbers $P_n$ existing in the series $\lambda(N)$ and it's not one of the multiples of any prime numbers less than $P_n$.

If $\lambda(N)$ is the summation till the number $N$, I will assume that the number $N$ is a very large prime numbers integer, whereas:

$$N = P_{n+1}, \quad P_{n+1} \to \infty$$

Assuming also that $N$ is very large so that it achieves the relation no. 22, thus we can get the following equation:

$$\lambda(P_{n+1}) = \frac{P_{n+1}^{1-s}}{1-s} \quad \ldots\ldots\ldots (27)$$

Also whereas the summation ends at the prime numbers $P_{n+1}$:

$$\lambda_{P_{n+1}}(P_{n+1}) = P_{n+1}^{-s} \quad \ldots\ldots\ldots (28)$$

But the series $\lambda_{Pn}$ remains and does not increase than the number $P_n$, because the numbers that follows $P_n$ is the square of $P_n$, as it is obvious through equation no.26, which is greater than $P_{n+1}$. So, the series $\lambda_{Pn}$ will be in the following form:

$$\lambda_{P_n}(P_{n+1}) = P_n^{-s} \quad \ldots\ldots\ldots (29)$$

Now, through the relation no. 23:

$$\frac{\lambda_{P_n}(P_{n+1})}{\lambda(P_{n+1})} = \frac{1}{P_n}\left(1-\frac{1}{2}\right)\left(1-\frac{1}{3}\right)\left(1-\frac{1}{5}\right)\cdots\left(1-\frac{1}{P_{n-1}}\right) \quad \ldots\ldots\ldots (30)$$

Also:

$$\frac{\lambda_{P_{n+1}}(P_{n+1})}{\lambda(P_{n+1})} = \frac{1}{P_{n+1}}\left(1-\frac{1}{2}\right)\left(1-\frac{1}{3}\right)\left(1-\frac{1}{5}\right)\cdots\left(1-\frac{1}{P_n}\right) \quad \ldots\ldots\ldots (31)$$

By substituting with the left hand side in the last two equations 30 and 31 of equations 28 and 29, we get the next two relations:

$$\frac{P_{n+1}^{-s}}{\frac{1}{1-s}P_{n+1}^{1-s}} = \frac{1}{P_{n+1}}\left(1-\frac{1}{2}\right)\left(1-\frac{1}{3}\right)\left(1-\frac{1}{5}\right)\cdots\left(1-\frac{1}{P_n}\right) \quad \ldots\ldots\ldots (32)$$

$$\frac{P_n^{-s}}{\frac{1}{1-s}P_{n+1}^{1-s}} = \frac{1}{P_n}\left(1-\frac{1}{2}\right)\left(1-\frac{1}{3}\right)\left(1-\frac{1}{5}\right)\cdots\left(1-\frac{1}{P_{n-1}}\right) \quad \ldots\ldots\ldots (33)$$

Assuming that the difference between the last two prime numbers $P_{n+1}$ and $P_n$ is $\Delta P_n$, and it can be described through the following equation:

$$P_{n+1} = P_n + \Delta P_n \quad \ldots\ldots\ldots (34)$$

By substituting with $P_{n+1}$ in the numerator of the left hand side in equation 32 from equation 34, we get the following relation:

$$\frac{(1-s)(P_n + \Delta P_n)^{-s}}{P_{n+1}^{1-s}} = \frac{1}{P_n + \Delta P_n}\left(1-\frac{1}{2}\right)\left(1-\frac{1}{3}\right)\cdots\left(1-\frac{1}{P_n}\right) \quad \ldots\ldots (35)$$

By subtracting equation 33 from the equation 35 we get the following:

$$(1-s)\left\{\frac{(P_n+\Delta P_n)^{-s}}{P_{n+1}^{1-s}} - \frac{P_n^{-s}}{P_{n+1}^{1-s}}\right\} = \left(1-\frac{1}{2}\right)\left(1-\frac{1}{3}\right)\cdots\left(1-\frac{1}{P_{n-1}}\right)\left\{\frac{1}{P_n+\Delta P_n}\left(1-\frac{1}{P_n}\right) - \frac{1}{P_n}\right\}$$

$$\Rightarrow (1-s)\frac{P_n^{-s}\left(1+\frac{\Delta P_n}{P_n}\right)^{-s} - P_n^{-s}}{P_{n+1}^{1-s}} = \left(1-\frac{1}{2}\right)\left(1-\frac{1}{3}\right)\ldots\left(1-\frac{1}{P_{n-1}}\right)\frac{1}{P_n}\left\{\frac{1}{1+\frac{\Delta P_n}{P_n}}\left(1-\frac{1}{P_n}\right) - 1\right\} \ldots (36)$$

By using Taylor series, the next term can be expand as follows:

$$\left(1+\frac{\Delta P_n}{P_n}\right)^{-s} = 1 - s\frac{\Delta P_n}{P_n} + \frac{1}{2}s(1+s)\left(\frac{\Delta P_n}{P_n}\right)^2 - \cdots$$

$$\Rightarrow P_n^{-s}\left\{\left(1+\frac{\Delta P_n}{P_n}\right)^{-s} - 1\right\} = P_n^{-s}\left\{\left(1 - s\frac{\Delta P_n}{P_n} + \frac{1}{2}s(1+s)\left(\frac{\Delta P_n}{P_n}\right)^2 - \cdots\right) - 1\right\}$$

$$\Rightarrow P_n^{-s}\left\{\left(1+\frac{\Delta P_n}{P_n}\right)^{-s} - 1\right\} = P_n^{-s}\left\{-s\frac{\Delta P_n}{P_n} + \frac{1}{2}s(1+s)\left(\frac{\Delta P_n}{P_n}\right)^2 - \cdots\right\}$$

Since the difference between any two very large respectively prime numbers must be much less than any of them, so the $\Delta P_n$ is much less than $P_n$, we can conclude the following limit:

$$\lim_{P_n \to \infty}\left\{-s\frac{\Delta P_n}{P_n} + \frac{1}{2}s(1+s)\left(\frac{\Delta P_n}{P_n}\right)^2 - \cdots\right\} = -s\frac{\Delta P_n}{P_n} \ldots\ldots\ldots (37)$$

Also the next part can be expand as follows:

$$\frac{1}{1+\frac{\Delta P_n}{P_n}} = 1 - \frac{\Delta P_n}{P_n} + \left(\frac{\Delta P_n}{P_n}\right)^2 - \left(\frac{\Delta P_n}{P_n}\right)^3 + \cdots \qquad \ldots\ldots\ldots (38)$$

Where $P_n$ is a very large number that approaches to infinity:

$$1 - \frac{1}{P_n} \to 1$$

Accordingly, the final bracket in equation 36 can be shortened as follows:

$$\frac{1}{1+\frac{\Delta P_n}{P_n}}\left(1-\frac{1}{P_n}\right) - 1 = -\frac{\Delta P_n}{P_n} + \left(\frac{\Delta P_n}{P_n}\right)^2 - \left(\frac{\Delta P_n}{P_n}\right)^3 + \cdots$$

When $P_n$ approaches to infinity, the last equation becomes as follows:

$$\frac{1}{1+\frac{\Delta P_n}{P_n}}\left(1-\frac{1}{P_n}\right) - 1 = -\frac{\Delta P_n}{P_n} \ldots\ldots\ldots (39)$$

Therefore, after using the relations 37, 38, and 39 by substitution in equation 36, it becomes as follows:

$$\frac{s(1-s)P_n^{-s}}{P_{n+1}^{1-s}} = \left(1-\frac{1}{2}\right)\left(1-\frac{1}{3}\right)\ldots\left(1-\frac{1}{P_{n-1}}\right)\frac{1}{P_n} \qquad \ldots\ldots\ldots (40)$$

By substituting $P_{n+1}$ which exists in the denominator of the left hand side of equation 40 in equation no. 34, results the following equations:

$$\frac{s(1-s)P_n^{-s}}{(P_n + \Delta P_n)^{1-s}} = \left(1-\frac{1}{2}\right)\left(1-\frac{1}{3}\right)\ldots\left(1-\frac{1}{P_{n-1}}\right)\frac{1}{P_n}$$

$$\Rightarrow \frac{s(1-s)P_n^{-s}}{P_n^{1-s}\left(1+\frac{\Delta P_n}{P_n}\right)^{1-s}} = \left(1-\frac{1}{2}\right)\left(1-\frac{1}{3}\right)\ldots\left(1-\frac{1}{P_{n-1}}\right)\frac{1}{P_n}$$

$$\Rightarrow \frac{s(1-s)}{\left(1+\frac{\Delta P_n}{P_n}\right)^{1-s}} = \left(1-\frac{1}{2}\right)\left(1-\frac{1}{3}\right)\ldots\left(1-\frac{1}{P_{n-1}}\right) \qquad \ldots\ldots\ldots (41)$$

Once again, by using the relations no. 37, 38, and 39 and substituting them in the last equation, results the following equation:

$$s(1-s)\left\{1-(1-s)\frac{\Delta P_n}{P_n}\right\} = \left(1-\frac{1}{2}\right)\left(1-\frac{1}{3}\right)\ldots\left(1-\frac{1}{P_{n-1}}\right) \quad \ldots\ldots\ldots(42)$$

Based on Riemann functional equation and as I previously declared, the non-trivial zeros are the zeros of the functions $\zeta(1-t)$ and $\zeta(t)$ together. Whereas they have the same function base, we can replace 1-s by s in equation 42, when s is one of the non-trivial zeros of zeta function; which makes both $\zeta(1-s)$ and $\zeta(s)$ equal to zero. Whereas $\zeta(1-s)$ achieve the same equations which we performed on $\zeta(s)$, we can replace 1-s by s and vice versa in equation 42 to come up with the following:

$$s(1-s)\left\{1-s\frac{\Delta P_n}{P_n}\right\} = \left(1-\frac{1}{2}\right)\left(1-\frac{1}{3}\right)\ldots\left(1-\frac{1}{P_{n-1}}\right) \quad \ldots\ldots\ldots(43)$$

By subtracting equation no. 42 from equation no. 43, results the following:

$$s(1-s)(-1+2s)\frac{\Delta P_n}{P_n} = 0$$

Since $\Delta P_n$ does not equal zero:

$$s(1-s)(-1+2s) = 0$$

Because the analysis was apply on eta functions, the initial or primary solution of this equation shows the real part of zeros of eta function:

$$Re[t] = 0, 1, 1/2$$

The zero is the real part of the classical zeros of $\eta_m(1-t)$ functions, while 1 is the real part of the classical zeros of $\eta_m(t)$ functions, and the half is the real part of the non-trivial zeros of all the functions forms of $\eta_m(1-t)$ and $\eta_m(t)$, which is the same non-trivial zeros of zeta function.

## 6) *Base of the imaginary part of the non-trivial zeros of zeta function:*

The next step is to perform a mathematical analysis on the functions of $\lambda$ after allocating N value by equalizing it with $P_{n+1}$, because the values of $\lambda$ function in mainly depends on the prime numbers, and we figured out equations that define its features. In this way, the equations that represent the imaginary part of zeros which we are going to have will be as afunction of the prime numbers. But the equations include all the natural numbers, so we'll generalize the solutions which depend on the prime numbers after obtaining all natural numbers with the same values which approach to those prime numbers.

Now and through equation no. 40:

$$\frac{s(1-s)P_n^{-s}}{P_{n+1}^{1-s}} = \left(1-\frac{1}{2}\right)\left(1-\frac{1}{3}\right)\ldots\left(1-\frac{1}{P_{n-1}}\right)\frac{1}{P_n}$$

With replacing 1-s by s and vice versa in the same equation, we get the following equation:

$$\frac{s(1-s)P_n^{s-1}}{P_{n+1}^{s}} = \left(1-\frac{1}{2}\right)\left(1-\frac{1}{3}\right)\ldots\left(1-\frac{1}{P_{n-1}}\right)\frac{1}{P_n} \quad \ldots\ldots\ldots(44)$$

By dividing equations no. 40 and 44, we get rid of the two factors s and (1-s) which cause the emergence of the classical zeros of eta function, and the factor causing the non-trivial zeros in the resulting equation will remains as follows:

$$\frac{P_{n+1}^s P_n^{-s}}{P_n^{s-1} P_{n+1}^{1-s}} = 1 \Rightarrow \frac{P_{n+1}^{2s-1}}{P_n^{2s-1}} - 1 = 0$$

$$\Rightarrow \frac{P_{n+1}^{s-1/2}}{P_n^{s-1/2}}\left\{\left(\frac{P_{n+1}}{P_n}\right)^{s-1/2} - \left(\frac{P_{n+1}}{P_n}\right)^{-(s-1/2)}\right\} = 0$$

$$\frac{1}{2i}\left\{e^{i(s-1/2)\ln\left(\frac{P_{n+1}}{P_n}\right)^{-i}} - e^{-i(s-1/2)\ln\left(\frac{P_{n+1}}{P_n}\right)^{-i}}\right\} = 0$$

$$\Rightarrow \sin\left\{(s-1/2)\ln\left(\frac{P_{n+1}}{P_n}\right)^{-i}\right\} = 0$$

$$\Rightarrow -i(s-1/2)\ln\left(\frac{P_{n+1}}{P_n}\right) = \pm k\pi, k = 1,2,3,\ldots$$

$$\Rightarrow s = \frac{1}{2} \pm i\frac{k\pi}{\ln\left(\frac{P_{n+1}}{P_n}\right)}, k = 1,2,3,\ldots \quad \ldots\ldots\ldots(45)$$

"This equation proves the Riemann hypothesis".

Now I will define $y(P_n)$ as being the imaginary part of the non-trivial zeros of zeta function, which takes the following form:

$$y(P_n) = \pm\frac{k\pi}{\ln\left(\frac{P_{n+1}}{P_n}\right)}, k = 1,2,3,\ldots$$

Through the substitution in equation no.34:

$$y(P_n) = \pm\frac{k\pi}{\ln\left(\frac{P_n + \Delta P_n}{P_n}\right)} = \pm\frac{k\pi}{\ln\left(1 + \frac{\Delta P_n}{P_n}\right)}$$

As for the natural logarithm in the denominator of the last equation, it can be expand by using Taylor Expansion as follows:

$$\ln\left(1 + \frac{\Delta P_n}{P_n}\right) = \frac{\Delta P_n}{P_n} - \frac{1}{2}\left(\frac{\Delta P_n}{P_n}\right)^2 + \frac{1}{3}\left(\frac{\Delta P_n}{P_n}\right)^3 - \ldots$$

Through the above, we can conclude the next limit when $\Delta P_n$ is much less than $P_n$:

$$\lim_{P_n \to \infty} \ln\left(1 + \frac{\Delta P_n}{P_n}\right) = \frac{\Delta P_n}{P_n} \quad \ldots\ldots\ldots(46)$$

$$y(P_n) = \pm\frac{k\pi P_n}{\Delta P_n} \quad \ldots\ldots\ldots(47)$$

Therefore:

If the function $\Pi(N)$ is Prime-counting function, which represents the number of prime numbers existing starting from 1 up to the number $N$, this function at large numbers approximates the following form:

$$\Pi(N) = \frac{N}{\ln N} \quad \ldots\ldots\ldots(48)$$

When $N$ is the prime number $P_n$, $\Pi(N)$ equals the order of prime number, i.e. it equals $n$, We can be prescribed through the following relation:

$$n = \Pi(P_n) \quad \ldots\ldots\ldots(49)$$

Which mean that $\Pi(N)$ is the inverse function of $P_n$, and it can be prescribed through the following relation:

$$P_n = \Pi^{-1}(n) \quad \ldots\ldots\ldots(50)$$

Since $\Pi^{-1}(n)$ is the inverse function of $\Pi(N)$:

$$\Pi(P_n) = \frac{P_n}{\ln P_n}$$

Thus:

$$n = \frac{P_n}{ln\, P_n} \quad \ldots\ldots\ldots (51)$$

Now I will derive the equation 51 through the next steps:

$$n + 1 = \frac{P_{n+1}}{ln\, P_{n+1}}$$

By substitution in equation no. 34:

$$n + 1 = \frac{P_n + \Delta P_n}{ln(P_n + \Delta P_n)} = \frac{P_n + \Delta P_n}{ln\left\{P_n\left(1 + \frac{\Delta P_n}{P_n}\right)\right\}} = \frac{P_n + \Delta P_n}{ln\, P_n + ln\left(1 + \frac{\Delta P_n}{P_n}\right)}$$

By substitution the logarithm in the equation 46:

$$n + 1 = \frac{P_n + \Delta P_n}{ln\, P_n + \frac{\Delta P_n}{P_n}} \quad \ldots\ldots\ldots (52)$$

By subtracting equation no. 51 from equation no. 52 we get the following:

$$1 = \frac{P_n + \Delta P_n}{ln\, P_n + \frac{\Delta P_n}{P_n}} - \frac{P_n}{ln\, P_n} = \frac{ln\, P_n\, (P_n + \Delta P_n)}{ln\, P_n\left(ln\, P_n + \frac{\Delta P_n}{P_n}\right)} - \frac{P_n\left(ln\, P_n + \frac{\Delta P_n}{P_n}\right)}{ln\, P_n\left(ln\, P_n + \frac{\Delta P_n}{P_n}\right)}$$

After removing the similar parts in the right hand side, we obtain the following equation:

$$\frac{\Delta P_n\, (ln\, P_n - 1)}{(ln\, P_n)^2\left(1 + \frac{\Delta P_n}{P_n\, ln\, P_n}\right)} = 1 \quad \ldots\ldots\ldots (53)$$

Since $P_n$ is a very large number:

$$ln\, P_n - 1 \to ln\, P_n$$

Since $\Delta P_n$ is much less than $P_n$:

$$1 + \frac{\Delta P_n}{P_n\, ln\, P_n} \to 1$$

After taking the previous two limits into consideration, the left hand side of equation no.53 will lead to the following form:

$$\frac{\Delta P_n\, (ln\, P_n - 1)}{(ln\, P_n)^2\left(1 + \frac{\Delta P_n}{P_n\, ln\, P_n}\right)} \to \frac{\Delta P_n}{ln\, P_n}$$

Thus, the equation 53 will lead to the following form:

$$\frac{\Delta P_n}{ln\, P_n} \to 1 \Rightarrow \Delta P_n = ln\, P_n \quad \ldots\ldots\ldots (54)$$

By dividing the tow hand sides of equation into $P_n$ results the following equation:

$$\frac{\Delta P_n}{P_n} = \frac{ln\, P_n}{P_n}$$

Through substituting the equation 51 in the right hand side of the last equation, we get the following result:

$$\frac{\Delta P_n}{P_n} = \frac{1}{n} \quad \ldots\ldots\ldots (55)$$

Upon this result, we can write the equation no. 47 in the following form:

$$y(P_n) = \pm nk\pi \quad \ldots\ldots\ldots (56)$$

The next step is to generalize equation no.56 to include all zeros that are being formed by all integers which approximate the infinity. Through the substituting with $P_n$ in equation 56 from equation 50, we get the following form:

$$y(\Pi^{-1}(n)) = \pm nk\pi \quad \ldots\ldots\ldots(57)$$

Whereas:

$$\Pi^{-1}(\Pi(n)) = n$$

Now, by substitution $\Pi(n)$ in place $n$ in the tow hand sides of equation 57, we get the following form:

$$y(n) = \pm k\pi \Pi(n)$$

Now, by the substituting in equation no.53, results the general formula of the imaginary part of non-trivial zeros at numbers that approach to infinity, as follows:

$$y(n) = \pm k\pi \frac{n}{\ln n}, k = 1,2,3,\ldots \ldots\ldots\ldots(58)$$

Since the number $k$ is an arbitrary number, this number must be defined by a certain value. The remained numbers are a repetition for the same value of imaginary part of non-trivial zeros at the number $n$.

❖ *Identifying the value of k:*

Through examining the equation 56, we can notice that the base of the imaginary part of non-trivial zeros at the prime numbers depends on two variables:

1) The number $n$ which represents the order of prime number.

2) The number $k$.

Where $k$ is an arbitrary number, while $n$ is an original variable, the number $k$ must be defined by a certain value. Although the other numbers that $k$ has, however it achieve the solution of the equation, they are just considered as repetition for the solution at the number $n$.

When examining the equation no. 58, we can notice that the base of the imaginary part of non-trivial zeros of zeta function at the large numbers which approaches to infinity leads to a similar form for the base of the classical zeros of eta function.

If the general base of the classical zeros of functions $\eta_m(1\text{-}t)$ and $\eta_m(t)$ is given through the relation no.15, its imaginary part is given through the following relation:

$$y(m) = \pm \frac{2j\pi}{\ln m}, \quad j = 1,2,3,\ldots \infty$$

Whereas $y(m)$ refers to the imaginary part of those zeros and $j$ here is an arbitrary number; while m in the original number.

The imaginary part of non-trivial zeros of zeta function that approaches to infinity is given through equation no. 58, so that the number $k$ is an arbitrary number and $n$ is the original number.

Since the base of the imaginary part of non-trivial zeros of zeta function at the large numbers approaches to infinity leads to the pass of the imaginary part of classical zeros of eta function, the two bases match each other when the following condition are achieved:

$$m = j = n$$

The other condition is:

$$k = 2$$

Accordingly, the path that the imaginary part of non-trivial zeros of zeta function leads to, at a very large numbers arises, and this path is determined through the following relation:

$$y(n) = \pm \frac{2\pi n}{\ln n} \quad \ldots\ldots\ldots (59)$$

Through the equation 59, the base of the non-trivial zeros of zeta function at numbers which approach the infinity will take the following form:

$$s = \frac{1}{2} \pm \frac{2\pi n}{\ln n} i \quad \ldots\ldots\ldots (60)$$

This result which proved by this analyses, matching results of some analyses that have been previously published. [10]

7) *Identifying the locations where the prime numbers exist through the non-trivial zeros:*

It was clearly indicated by the recent analyses that the imaginary part of the non-trivial zeros of zeta function does not always refer to the locations of the prime numbers, as every non-trivial zero does not relate to a certain order of prime numbers, but when the condition stating that the $y(n)$ value must approach to infinity is achieved, the prime numbers are located where the ratio of $y(n)$ to $2\pi$ is an integer, and this integer is equal to $\Pi(n)$, and as shown prescribed in equation 49, the number $n$ in equation no. 59 (which achieves this condition) is order of prime number.

Finally, the relations which we have reached become totally valid and not just approximately relations at a very large numbers that approach to infinity, but although they become less accurate when the numbers on which they depend becomes smaller, However, it gives a vision for path of non-trivial zeros of zeta function.

Accordingly, the zeros at very large numbers take a relation as afunction of the Prime-counting function $\Pi(n)$ as follows:

$$s = \frac{1}{2} \pm 2\pi\Pi(n)i \quad \ldots\ldots\ldots (61)$$

## CONCLUTION:

This analysis has proven the Riemann Hypothesis, since the real part of the non-trivial zeros of zeta function is always constant and equals to half.

This analysis has proven also that the imaginary part of the non-trivial zeros of zeta function is somehow connected to the Prime-counting function $\Pi(n)$, and it is directly proportional with $\Pi(n)$ at $n$ approaches to infinity.

This analysis has proven as well that the imaginary part of the non-trivial zeros of zeta function refers somehow to the locations where the prime numbers exist according to the condition achieving the equation no. 49.